\documentclass[a4paper, 9pt]{article}

\usepackage[applemac]{inputenc}
\usepackage{bibtopic}
\usepackage{amsmath,amssymb,amsthm}

\renewcommand{\>}{\rangle}
\newcommand{\<}{\langle}

\newtheorem{thm}{Theorem}[section]

\newtheorem{lem}{Lemma}[section]

\title{{A remark on the definability of the Fitting subgroup and the soluble radical}}
\date{}
\author{A. Ould Houcine\thanks{The author is supported by SFB 878.}}
\begin{document}

\maketitle

\begin{abstract} Let $G$ be an arbitrary  group. We show that if the Fitting subgroup of $G$  is  nilpotent then it is definable.  We show also that the  class of groups whose  Fitting subgroup is nilpotent of class at most $n$  is elementary.  We give an example of a  group (arbitrary saturated) whose  Fitting subgroup  is definable but not nilpotent. Similar results for the soluble radical are  given. 

\end{abstract}

\section{Introduction}

Let $G$ be a group.  The \textit{Fitting subgroup} of $G$ is the subgroup generated by all normal nilpotent subgroups of $G$. It  will be denoted  $F(G)$.   The subgroup generated by all the normal soluble subgroups of $G$ is called the \textit{soluble radical} and will be denoted  $R(G)$. 

We see that if $G$ is finite then $F(G)$ is nilpotent and $R(G)$ is soluble but in general $F(G)$ is not necessarily nilpotent, similarly $R(G)$ is not necessarily soluble. However, we can see that for an arbitrary group $G$, $F(G)$ is  locally nilpotent and $R(G)$ is  locally soluble. 

It is well-known that the Fitting subgroup  (or the soluble radical) is definable in many classes of groups.  The definability of the Fitting subroup of a group $G$  makes it   more accessible for using model-theoretic techniques. 

The purpose of the present note is to give  relations between the nilpotency and the definability of the Fitting subgroup (we treat also the soluble radical). The first observation is that the nilpotency (resp.  the solubility) of the Fitting subgroup (resp. of the soluble radical) implies its definability. 

\begin{thm} \label{mainthm1}$\;$

$(1)$ For any $n\geq 1$, there exists a formula $\phi_n(x)$, depending only on $n$, such that for any group $G$,  if $F(G)$ is nilpotent of class $n$ then it is definable by $\phi_n$. 

$(2)$ For any $n\geq 1$, there exists a formula $\psi_n(x)$, depending only on $n$, such that for any group $G$,  if $R(G)$ is of derived length  $n$ then it is definable by $\psi_n$. 
\end{thm}

One may ask  if the given  formula can define the Fitting subgroup in elementary extensions of $G$.  We show that it is  the case. 

\begin{thm}  \label{mainthm2} The  class of groups whose  Fitting subgroup is nilpotent of class at most $n$  is elementary.  Similarly for the class of groups whose soluble radical is of derived length at most $n$. 
\end{thm}

One of the consequences of Theorem \ref{mainthm1} and \ref{mainthm2} is that if $G$ is an $\mathfrak M_C$-group then for any group $K$ elementary equivalent to $G$, $F(K)$ is nilpotent and definable. Indeed by a result of Derakhshan and Wagner \cite{Der-Wag, wag} $F(G)$ is nilpotent and we apply Theorem \ref{mainthm1} and \ref{mainthm2}. We note that being an $\mathfrak M_C$-group is not necessarily conserved by elementary equivalence.  Recently Altıinel and Baginski \cite{al-bag}  showed that in an $\mathfrak M_C$-group every nilpotent subgroup is contained in a definable nilpotent subgroup. 

As noticed above $F(G)$ is locally nilpotent. In fact $F(G)$ satisfies  a stronger property: the normal closure (in $G$)  of any finite subset of $F(G)$ is nilpotent. To formulate the next result we introduce the following definition: we say that $F(G)$ is \textit{uniformly normaly locally nilpotent} if for any $m$ there exists $d(m)$ such that for any finite subset $A \subseteq F(G)$ of cardinality at most $m$  the normal closure of $A$ is nilpotent of class at most $d(m)$.  The notion of $R(G)$ is  \textit{uniformly normaly locally soluble} is defined in a similar way. 

\begin{thm} \label{mainthm3} Let $G$ be an $\aleph_0$-saturated group. Then $F(G)$ is definable if and only if $F(G)$ is uniformly normaly locally nilpotent. Similarly $R(G)$ is definable if and only if $R(G)$ is uniformly normaly locally soluble. 
\end{thm}

One may also ask if the converse of Theorem \ref{mainthm1} is true, that is   if the definability of the Fitting subgroup implies its nilpotency (similarly for the soluble radical). At the end of the next section we give an example, based on a result of S. Bachmuth and H.Y. Mochizuki \cite{Bach-Moch}, of a group $G$ (even $\aleph_0$-saturated) such that $F(G)$ is definable but not nilpotent,  similarly for the soluble radical.

\bigskip
\noindent \textbf{Acknowledgement.} I am  very grateful  to Tuna Altıinel for  stimulating discussions during the conference "Model Theory of Groups
(Interaction between Model Theory and Geometric Group Theory)"
CIRM, Luminy, November 21-25, 2011.

\section{Proofs}

The notations which we are going to use are very standard. Given a group $G$ and $g, h \in G$ the commutator $[g,h]$ is defined to be $g^{-1}h^{-1}gh$.  We use the notation $g^h=h^{-1}gh$. If $A$ is a subset of $G$, then $\<A\>$ denotes the subgroup generated by $A$ and $\<A\>^G$ stands for  the normal closure of $A$. For $g \in G$, $g^G$ denotes the conjugacy class of $g$. We see that $\<g^G\>=\<g\>^G$.  If $H, K$ are subgroups of $G$ then $[H,K]$ is the subgroup generated by the commutators $[h,k]$ where $h \in H, k \in K$.

The next lemma can be deduced from \cite{rob}  but for the  sake of completeness we provide  a proof. 

\begin{lem} \label{lem1} Let $G$ be a group,  $H$ and $K$ be normal subgroups of $G$. Suppose that $H$ is generated by $A$ and $K$ is generated by $B$. Then  $[H,K]$ is generated by $X=\{[a^\alpha,b^\beta] : a \in A, b\in B, \alpha, \beta \in G\}$. 
\end{lem}

\proof We note that $\<X\>$ is normal and $\<X\> \leq [H,K]$. We must show that for any  $h \in H, k \in K$, $[h,k] \in \<X\>$. There exists a word  $w$ (resp. $v$) over the alphabet  $A^{\pm 1}$ (resp. over $B^{\pm 1}$) such that  $h=w(\bar a), k=v(\bar b)$.  The result follows by  induction on  $|w|+|v|$, where $|.|$  denotes the word-length, by using the following formulae:
$$
[x,yz]=[x,z]  [x,y]^z, \; \; [xy,z]=[x,z]^y[y,z], \;\; [x,y]^z=[x^z, y^z], $$$$ [x^{-1}, y]=[x,y^{x^{-1}}]^{-1}, \;\; [x,y^{-1}]=[x^{y^{-1}}, y]^{-1}. \qed
$$

\bigskip

We define  by induction on $n$, the word  $u_{n}(x_1, \dots, x_{n+1})$ as follows:  
$$
u_{1}(x_1, x_2)=[x_1, x_2], 
$$
$$
u_{n+1}(x_1, \dots, x_{n+1}, x_{n+2})=[u_{n}(x_1, \dots, x_{n+1}), x_{n+2}].
$$

We see that a group $G$ is nilpotent of class at most $n$ if and only if $G \models \forall x_1  \dots \forall x_{n+1} u_n(x_1, \dots, x_{n+1})=1$. 

Similarly,  by induction on $n$, define the word $v_{n}(x_1, \dots, x_{2^n})$   as follows:
$$
v_{1}(x_1, x_2)=[x_1, x_2], v_{2}(x_1, x_2, x_3, x_4)=[[x_1,x_2],[x_3, x_4]], 
$$
$$
v_{n+1}(x_1,  \dots, x_{2^n}, x_{1+2^n}, \dots, x_{2^{n+1}})=[v_n(x_1, \dots, x_{2^n}), v_n(x_{1+2^n}, \dots, x_{2^n+2^n})]. 
$$

As above  a group $G$ is soluble of derived length  at most $n$ if and only if $G \models \forall x_1  \dots \forall x_{2^n} v_n(x_1, \dots, x_{2^n})=1$. 

\bigskip
\begin{lem} \label{lem2} Let $N$ be a normal subgroup of $G$ generated by  $B$. 

$(1)$ The $n$-th term  $N^n$ of the descending central series of $G$ is generated  by 
$$
\{u_{n}(b_1^{\alpha_1},\dots, b_{n+1}^{\alpha_{n+1}} ):  b_i\in B, \alpha_i \in G\}. 
$$

 $(2)$ The $n$-th term  $N^{(n)}$ of the derived series of $G$ is generated by 
 
 $$
\{v_{n}(b_1^{\alpha_1},\dots, b_{2^n}^{\alpha_{2^n}} ):  b_i\in B, \alpha_i \in G\}. 
$$

\end{lem}

\proof  We treat only $(1)$ the proof of case $(2)$ proceeds in a similar way. The proof is by induction on $n$. For  $n=1$, we get $N^1=N'$ and the result is a consequence of Lemma \ref{lem1}. For  $n+1$,  $N^{n+1}=[N^n, N]$.  By induction $N^n$ is generated by 
$$
\{u_{n}(b_1^{\alpha_1},\dots, b_{n+1}^{\alpha_{n+1}} ):  b_i\in B, \alpha_i \in G\}.
$$

By Lemma \ref{lem1}, since  $N^n$  is normal in $G$,  $N^{n+1}$  is generated by  

$$
\{[u_{n}(b_1^{\alpha_1},\dots, b_{n+1}^{\alpha_{n+1}})^{\beta}, b_{n+2}^{\alpha_{n+2}}] :  b_i\in B,  \beta, \alpha_i \in G\}, 
$$
which  coincides with the following  set
$$
\{u_{n+1}(b_1^{\alpha_1},\dots, b_{n+1}^{\alpha_{n+1}}, b_{n+2}^{\alpha_{n+2}}) :  b_i\in B, \alpha_i \in G\}. \qed
$$

\begin{lem} \label{lem3} For any $(n,m) \in \mathbb N^2$, there exists a formula $\varphi_{n,m}(x_1, \dots, x_m)$ satisfying  the following property:  for any group $G$, for any finite subset $\{b_1, \dots, b_m\}$ of  $G$, $G \models \varphi_{n,m}(b_1, \dots, b_m)$ if and only if  $\<\{b_1, \dots, b_m\}\>^G$ is nilpotent of class at most $n$.  Similarly for  the soluble case. 
\end{lem}

\proof  We treat the nilpotent case,  the proof of the soluble case  proceeds in a similar way. Let us denote by $S_{nm}$ the set of all applications from the set $\{1, \dots, n+1\}$ to the set $\{1, \dots, m\}$.  Let $\varphi_{n,m}(x_1, \dots, x_{m})$ to  be the following formula
$$
\forall y_1 \forall y_2 \dots \forall y_{n+1}\bigwedge_{ \sigma \in S_{nm}}u_{n}(x_{\sigma(1)}^{y_1}, \dots, x_{\sigma(n+1)}^{y_{n+1}})=1.
$$

Let $A=\{b_1, \dots, b_m\} \subseteq G$. Then $\<A\>^G=\<\{b_i^{g} | 1 \leq i\leq m,  g \in G\}\>$.  By Lemma \ref{lem2}  $(\<A\>^G)^n$ is generated by 
$$
\{u_{n}(b_{\sigma(1)}^{g_1}, \dots, b_{\sigma(n+1)}^{g_{n+1}}), \sigma \in S_{nm},  g_i \in G\},
$$
which gives the required result.  \qed

\bigskip
\noindent \textbf{Proof of  Theorem \ref{mainthm1}.} We treat only $(1)$ the proof of case $(2)$ proceeds in a similar way. Let $n$ to be the nilpotency class of $F(G)$. We have  $a \in F(G)$ if and only if   $\<a\>^G$ is nilpotent of class at most  $n$. But $\<a\>^G$ is generated by  $\{a^g : g \in G\}$ and thus we conclude by Lemma \ref{lem2} that  
$$
a \in F(G) \hbox{ if and only if  } G \models \forall x_1  \cdots \forall x_{n+1}( u_{n}(a^{x_1}, \dots, a^{x_{n+1}})=1). 
\qed$$

\bigskip
\noindent \textbf{Proof of Theorem \ref{mainthm2}.} We treat only the nilpotent case,   the proof of the soluble case  proceeds in a similar way.  Let $p \geq 0$ and 
$$
\phi_{n,m}=: \forall x_1 \dots \forall x_m (\bigwedge_{1 \leq i \leq m}\varphi_{n,1}(x_i) \implies \varphi_{p,m}(x_1, \dots, x_m)),
$$
where $\varphi_{n,m}$ is the formula given by Lemma \ref{lem3}. Let 
$$
T_p=\{\phi_{n,m}|n,m \in \mathbb N^*\}.
$$
 
 We claim that  $G \models T_p$ if and only if  $F(G)$ is nilpotent of  class at most $p$. Let $G \models T_p$.  We must  show that for any $b_1, \cdots, b_{p+1} \in F(G)$, $F(G) \models  u_p(b_1, \dots, b_{p+1})=1$. 
 
 Let $b_1, \cdots, b_{p+1} \in F(G)$.  For each $1 \leq i \leq p+1$, since $b_i \in F(G)$, there exists $n_i$ such that $\<b_i\>^G$ is nilpotent of class $n_i$.  Let $n=\max_{1 \leq i \leq {p+1}}n_i$.  Hence, we have $G\models \bigwedge_{1 \leq i \leq p+1}\varphi_{n,1}(b_i)$.  Therefore, $G \models \varphi_{p,p+1}(b_1, \dots, b_{p+1})$ and thus $\<b_1, \dots, b_{p+1}\>^G$ is nilpotent of class at most $p$ which gives the required result.

 Now, we show that if $F(G)$ is nilpotent of class at most $p$, then $G \models T_p$.  Let $n, m \in \mathbb N^*$ and let $b_1, \dots, b_m$ such that $G \models \bigwedge_{1 \leq i \leq m}\varphi_{n,1}(b_i)$.  Hence, for each $i$, $\<b_i\>^G$ is nilpotent of class at most $n$. Hence $\<b_i\>^G \leq F(G)$  and thus $\<b_1, \dots, b_m\>^G$ is nilpotent of class at most $p$ and thus $G \models \varphi_{p,m}(b_1, \dots, b_m)$ as required.  \qed

 \bigskip $\;$
 
\noindent \textbf{Proof of Theorem \ref{mainthm3}.} Suppose that $G$ is $\aleph_0$-saturated. We will deal only with the Fitting subgroup, the case of the soluble radical is similar. Suppose that $F(G)$ is definable. Since for any $b_1, \dots, b_m \in F(G)$, $\<b_1, \dots, b_m\>^G$ is nilpotent and $G$ is $\aleph_0$-saturated we conclude, using Lemma \ref{lem3}, that there is an uniform bound on the nilpotency class of $\<b_1, \dots, b_m\>^G$.  Conversely, if there is an uniform bound on the nilpotency class of nilpotent subgroups of the form $\<b\>^G$, then $F(G)$ is definable by the formula given in the proof of Theorem \ref{mainthm1}. \qed

\bigskip
We end this section by giving an example showing that the definability of the Fitting subgroup does not imply necessarily its nilpotency, similarly for the soluble radical. 

Recall that the (right) $n$-Engel word $[x,_ny]$ is defined inductively by $[x,_1y]=[x,y]$ and $[x,_{n+1}y]=[[x,_ny], y]$.  A group $G$ is said to be $n$-Engel if $G \models \forall x \forall y [x,_ny]=1$.  It is known that a group $G$ is $2$-Engel if and only if the normal closure of any element of $G$ is abelian.  L. C. Kappe and W. P. 
Kappe \cite{Kap} proved that a group is 3-Engel if and only if the normal 
closure of any  element of $G$ is nilpotent of class at most 2. It follows that if $G$ is $3$-Engel then $F(G)=R(G)=G$.

S. Bachmuth and H.Y. Mochizuki \cite{Bach-Moch} constructed a nonsoluble $3$-Engel group $G$. We see that if $L$ is an elementary extension of $G$, which can be shosen arbitrary saturated,  then $L$ is also $3$-Engel, nonsoluble and $F(L)=R(L)=L$. Hence $L$ is an example of a (arbitrary saturated) goup in which $F(G)$ is definable (by the obvious formula $x=x$) but  not nilpotent, similarly for the soluble radical.  

If we  want to have an example in which the Fitting subgroup is a proper subgroup, we  can modify the previous example slightly as follows. We let $K=G \times S$ where $S$ is a finite nonabelian simple group. Then by picking a nontrivial element $s \in S$ we see that $G$ is the centralizer $C_K(s)$. We see also that $F(K)=R(K)=G$. Hence the Fitting subgroup of $K$ is a proper definable subgroup which is not nilpotent and even nonsoluble and the same property holds also in elementary extensions of $K$.

\bigskip
\small \noindent Abderezak OULD HOUCINE, \\
Mathematisches Institut und 
Institut für Mathematische Logik und Grundlagenforschung, 
Fachbereich Mathematik und Informatik,  
Universität Münster,  
Einsteinstrasse 62,  
48149 Münster,  
Germany. 

\noindent Universit\'e  de Lyon; Universit\'e Lyon 1; INSA de Lyon, F-69621; Ecole Centrale
de Lyon; CNRS, UMR5208, Institut Camille Jordan, 43 blvd du 11 novembre
1918, F-69622 Villeurbanne-Cedex, France.  \\

\noindent \textit{E-mail}:\textrm{ould@math.univ-lyon1.fr}

\end{document}